\magnification=\magstep1
\vbadness=10000
\hbadness=10000
\tolerance=10000

\proclaim What is the monster.

Richard E. Borcherds, 
\footnote{$^*$}{ Supported by 
NSF grant DMS-9970611.}

Mathematics department,
Evans Hall \#3840,
University of California at Berkeley,
CA 94720-3840
U. S. A.

e-mail: reb@math.berkeley.edu

www home page  www.math.berkeley.edu/\hbox{\~{}}reb

\bigskip

When I was a graduate student, my supervisor John Conway would bring
his 1 year old son into the department, who was soon known as the
baby monster.  A more serious answer to the question is that the
monster is the largest of the (known\footnote{**}{The
announcement of the classification of the finite simple groups about
20 years ago was a little over-enthusiastic, but a recent 
1300 page preprint by Aschbacher
and Smith should finally complete it.}) 
sporadic simple groups. Its name
comes from its size: the number of elements is 
$$\eqalign{
&8080,17424,79451,28758,86459,90496,17107,57005,75436,80000,00000 \cr
= 
&2^{46}.3^{20}.5^9.7^6.11^2. 13^3.17.19.23.29.31.41.47.59.71,\cr} $$
 about equal to the
number of elementary particles in the planet Jupiter.  

The monster was originally predicted to exist by Fischer and by Griess
in the early 1970's.  Griess constructed it a few years later 
in an extraordinary tour de force as the group of linear transformations
on a vector space of dimension 196883 that preserved a certain
commutative but non-associative 
bilinear product, now called the Griess product. 

Our knowledge of the structure and representations of the monster is
now pretty good. The 194 irreducible complex representations were
worked out by Fischer, Livingstone and Thorne (before the monster was
even shown to exist) They take up 8 large pages in the atlas [A] of
finite groups, which is the best single source of information about
the monster (and other finite simple groups). The subgroup structure
is mostly known; in particular there is an almost complete list of the
maximal subgroups, and the main gaps in our knowledge concern
embeddings of very small simple groups in the monster.  If anyone
wishes to multiply elements of the monster explicitly, R. A.  Wilson
can supply two matrices that generate the monster. But there is a
catch: each matrix takes up about 5 gigabytes of storage, and to quote
from Wilson's atlas page: ``standard generators have now been made as
$196882
\times 196882$ matrices over GF(2)... They have been multiplied
together, using most of the computing resources of Lehrstuhl D f\"ur
Mathematik, RWTH Aachen for about 45 hours...''.  (The difficulty of
multiplying two elements of the monster is caused not so much by its
huge size as by the lack of ``small'' representations; for
example, the symmetric group $S_{50}$ is quite a lot bigger than the
monster, but it only takes a few minutes to multiply two elements by
hand.)  Finally the modular representations of the monster for large
primes were worked out by Hiss and Lux; the ones for small primes
still seem to be out of reach at the moment.

In the late 1970's John McKay decided to switch from finite group theory
to Galois theory. One function that turns up in Galois theory is the
elliptic modular function
$$j(\tau)=q^{-1}+744+196884q+21493760q^2+\cdots=\sum c(n)q^n$$ ($q=e^{2\pi i
\tau}$), which is essentially the simplest non-constant function
invariant under the action $\tau\mapsto (a\tau+b)/(c\tau+d)$ of
$SL_2({\bf Z})$ on the upper half plane $\{\tau|\Im(\tau)>0\}$. He noticed
that the coefficient 196884 of $q^1$ was almost equal to the degree
196883 of the smallest complex representation of the monster (up to a
small experimental error). The term ``moonshine'' roughly means weird
relations between sporadic groups and modular functions (and anything
else) similar to this.  It was clear to many people that this was just
a meaningless coincidence; after all, if you have enough large
integers from various areas of mathematics then a few are going to be
close just by chance, and  John McKay  was told that his
observation was about as useful as looking at tea-leaves.  John
Thompson took McKay's observation further, and pointed out that the
next few coefficients of the elliptic modular function were also
simple linear combinations of dimensions of irreducible
representations of the monster; for example,
$21493760=21296876+196883+1$.  He  suggested that there
should be a natural infinite dimensional graded representation
$V=\sum_{n\in {\bf Z}}V_n$ of the monster such that the dimension of $V_n$
is the coefficient $c(n)$ of $q^n$ in $j(\tau)$, at least for $n\ne
0$. (The constant term of $j(\tau)$ is arbitrary as adding a constant
to $j$ still produces a function invariant under $SL_2({\bf Z})$, and
is set equal to 744 mainly for historical reasons.) 
Conway and Norton [C-N] followed up Thompson's
suggestion of looking at the McKay-Thompson series $T_g(\tau)=\sum_n
Trace(g|V_n)q^n$ whose coefficients are given by the traces of
elements $g$ of the monster on the representations $V_n$, and found by
calculating the first few terms that they all seemed to be Hauptmoduls
of genus 0. (A Hauptmodul is a  function similar to $j$, but invariant
under some group other than $SL_2({\bf Z})$).  Atkin, Fong, and Smith showed by computer calculation
that there was indeed an infinite dimensional graded representation of
the monster whose McKay-Thompson series were the
Hauptmoduls found by Conway and Norton, and 
soon afterwards  Frenkel, Lepowsky, and Meurman
explicitly constructed  this representation using vertex operators.

If a group acts on a vector space it is natural to ask if 
it preserves any algebraic structure, such as a bilinear form or product. 
The monster module constructed by FLM has a vertex algebra structure invariant under the action of the monster. Unfortunately there is no easy way 
to explain what a vertex algebra is; see [K] for the best introduction to them.
Vertex algebras are a generalization of
commutative rings with derivations (at least in characteristic 0). 
Roughly speaking they can be thought of as commutative rings with derivation
where the ring multiplication is not quite defined everywhere; 
this is analogous to rational maps in algebraic geometry, which are
also not quite defined everywhere. A more concrete but less intuitive
definition of a vertex algebra is that it consists of a space
with a countable number of bilinear products satisfying 
certain rather complicated identities. In the case of 
the monster vertex algebra $V=\oplus V_n$, this gives bilinear
maps from $V_i\times V_j$ to $V_k$ for all integers $i,j,k$, 
and the special case of the map from $V_2\times V_2$ to $V_2$ is
(essentially) the Griess product. So the Griess algebra is 
a sort of section of the monster vertex algebra. 

Following an idea of Frenkel, the monster vertex algebra can be used
to construct the monster Lie algebra by using the Goddard-Thorn
no-ghost theorem from string theory.  This is a ${\bf Z}^2$-graded Lie
algebra, whose piece of degree $(m,n)\in {\bf Z}^2$ has dimension
$c(mn)$ whenever $(m,n)\ne (0,0)$.  The monster should be thought of
as a group of ``diagram automorphisms'' of this Lie algebra, in the
same way that the symmetric group $S_3$ is a group of diagram
automorphisms of the Lie algebra $D_4$.  The monster Lie algebra has a
denominator formula, similar to the Weyl denominator formula for
finite dimensional Lie algebras and the Macdonald-Kac identities for
affine Lie algebras, which looks like
$$j(\sigma)-j(\tau)=p^{-1}\prod_{m>0\atop n\in
{\bf Z}}(1-p^mq^n)^{c(mn)}$$ where $p=e^{2\pi i \sigma}$, $q=e^{2\pi i
\tau}$. 
This formula was discovered independently in the 1980's by several people,
including Koike,
Norton, and Zagier. There are similar identities with $j(\tau)$ replaced
by the McKay-Thompson series of any element of the monster, and Cummins and Gannon
showed that any function satisfying such identities is a  Hauptmodul.
So this provides some sort of explanation of Conway and Norton's
observation that the McKay-Thompson series are all Hauptmoduls.  

So the question ``What is the monster'' now has several reasonable answers:
\item{1} The monster is the largest sporadic simple group, or alternatively
the unique simple group with its order. 
\item{2} It is the automorphism group of the Griess algebra. 
\item{3} It is the automorphism  group of the monster
vertex algebra. (This is probably the best answer.)
\item{4} It is a group of diagram automorphisms of
the monster Lie algebra.

Unfortunately none of these definitions is completely satisfactory.
At the moment all constructions of the algebraic structures above seem
 artificial; they are constructed as sums of two or more
apparently unrelated spaces, and it takes a lot of effort to define
the algebraic structure on the sum of these spaces and to check that
the monster acts on the resulting structure.

\proclaim References.

\item{[A]} 
J. H. Conway, R. T. Curtis, S. P. Norton,
   R. A. Parker, R. A. Wilson, Atlas of finite groups, Clarendon Press, Oxford,
   1985.
\item{[C-N]} 
J. H. Conway, S. Norton, Monstrous moonshine,
   Bull. London. Math. Soc. 11 (1979) 308-339.
\item{[K]}
V. Kac, Vertex algebras for beginners. Second edition. 
University Lecture Series, 10. American Mathematical Society, Providence, RI,
1998. 
\bye